\documentclass[12pt,twoside,doublespace]{article}

\usepackage{amssymb}
\usepackage{graphicx}

\def\smskip{\par\vskip 5 pt}
\def\QED{\hfill $\Box$\smskip}
\newtheorem{theorem}{Theorem}
\newtheorem{lemma}{Lemma}

\begin{document}

\begin{center}

\vspace{35pt}

{\Large \bf A Non-Monotone Conjugate Subgradient }

\vspace{5pt}

{\Large \bf  Type Method for Minimization of }

\vspace{5pt}

{\Large \bf  Convex Functions}

\vspace{5pt}

\vspace{35pt}

{\sc I.V.~Konnov\footnote{\normalsize E-mail: konn-igor@ya.ru}}

\vspace{35pt}

{\em  Department of System Analysis
and Information Technologies, \\ Kazan Federal University, ul.
Kremlevskaya, 18, Kazan 420008, Russia.}

\end{center}

\vspace{35pt}

\begin{abstract}

We suggest a conjugate subgradient type method without
any line-search for minimization of convex non differentiable
functions. Unlike the custom methods of this class,
it does not require monotone decrease of the goal function
and reduces the implementation cost of each iteration essentially.
At the same time, its step-size procedure takes into account behavior
of the method along the iteration points.
Preliminary results of computational experiments confirm efficiency of
the proposed modification.

\end{abstract}


\newpage

%1111111111111111111111111111111111111111111111111111111111111111111111111

\section{Introduction} \label{sc:1}

Let $f: \mathbb{R}^{n} \rightarrow \mathbb{R}$ be a convex and continuous, but not necessarily differentiable
function defined on the real $n$-dimensional Euclidean space $\mathbb{R}^{n}$. Then we can
consider the well-known problem of minimizing $f$ over $\mathbb{R}^{n}$, or briefly,
\begin{equation} \label{eq:1.1}
 \min _{x \in \mathbb{R}^{n}} \rightarrow f(x).
\end{equation}
We denote by $X^{*}$ and $f^{*}$ the set of solutions of problem (\ref{eq:1.1})
and the optimal value in (\ref{eq:1.1}), respectively. There exist a great number of
significant applications of convex minimization problems having just
non-differentiable goal functions; e.g. see \cite{Sho85,Pol87,MN92} and the references therein.
For this reason, their theory and methods were developed rather well. In particular,
many different iterative methods were proposed for finding solutions of
convex non-differentiable (non-smooth) minimization problems; e.g.
see \cite{Sho85,Kiw85,Pol87,MN92,HL93,Kon13} and the references therein.
We recall that most applications admit calculation of only one arbitrary taken element from the
subdifferential of $f$ at any point.

During rather long time, most efforts
were concentrated on developing more powerful and rapidly convergent methods within this setting,
such as space dilation and bundle type ones, which admit
complex transformations at each iteration.
However, significant areas of applications related to decision making in
industrial, transportation, information and communication systems,
having large dimensionality and inexact data together with
scattered necessary information force one to avoid
complex transformations and even line-search procedures and to
apply mostly simple methods, whose iteration computation expenses
and accuracy requirements are rather low. We observe that the problem of
creation of efficient low cost non-smooth optimization methods is more difficult
in comparison with that in the smooth case.

For instance, let us take the simplest subgradient method:
\begin{equation} \label{eq:1.2}
\begin{array}{c}
\displaystyle
 x^{k+1} :=x^{k} - \lambda _{k} g^{k}, g^{k} \in \partial f (x^{k}), \lambda _{k} > 0,
\end{array}
\end{equation}
where $\partial f(x)$ denotes the subdifferential of $f$ at $x$. The choice of
the step-size $\lambda _{k}$, which provides convergence to a solution, can
conform to various rules. The most popular is to take the so-called divergent series rule:
\begin{equation} \label{eq:1.3}
\begin{array}{c}
\displaystyle
  \lim \limits_{k\to \infty } \lambda _{k}=0, \quad \sum \limits_{k=0}^{\infty } \lambda _{k} = \infty.
\end{array}
\end{equation}
However, convergence of such a method is rather slow because of the second
condition in (\ref{eq:1.3}) and non-adaptive step-size choice.
There are several ways to speed up convergence of the subgradient methods
with the same storage and computation expenses per iteration.

One of them consists in calculating the step-size via
utilization of a priori information such as the optimal value $f^{*}$
or some condition numbers; e.g. see \cite{Sho85,Pol87}. Usually, one must
take only some their inexact estimates, which again leads to slow convergence.
The second way consists in creating simple descent or relaxation subgradient methods;
see \cite{Wol75,DV81,Che82,Kon82,Kiw83,Kon84a} and also \cite{Kiw85,MGN87,Kon01,Kon13}.
These methods demonstrate more stable convergence in comparison with the custom
subgradient methods and even attain a linear rate of convergence, but the requirement
of monotone decrease of the goal function
at each iteration leads to rather small step-size and may result in
low convergence.
The third way consists in utilization of averaging procedures
for subgradients from several previous
iteration points in order to find the current direction;
see \cite{Gup78,Che81,MGN87}. This property is also used in most
descent subgradient methods, but here no descent is required for
iteration points. As a result, one can obtain very flexible methods, which however
require mostly divergent series type rules for providing convergence, so that
their convergence rate estimate is usually the same as that of the custom subgradient method
(\ref{eq:1.2})--(\ref{eq:1.3}).

Therefore, one is interested in further development of these directions in order to
create more efficient versions of subgradient methods.
Recently,  a simple adaptive step-size procedure for smooth optimization methods
was proposed in \cite{Kon18,Kon18b,Kon18d}. In this paper, we follow this approach, but
due to the non-smooth goal function we combine it with some procedures and rules from
conjugate subgradient methods (see \cite{Lem75,Wol75}) and
non monotone averaging subgradient methods (see \cite[Chapter IV]{MGN87}).
Our method admits different changes of the step-size and wide variety
of implementation rules. It does not utilize
a priori information, but takes into account behavior of the iteration sequence.
Preliminary results of computational experiments confirm its efficiency.

%22222222222222222222222222222222222222222222222222222222222222222222222222222

\section{Basic Preliminaries}\label{sc:2}

We first recall two well-known auxiliary properties of vector sequences; see \cite{Lem75,Wol75} and
\cite[Section 4.3]{Kon01}.  Given  a set $X$,
we denote by  conv$X$  the convex hull of $X$ and
by  Nr$X$ the element of $X$ nearest to origin.
Also, we denote by $B(x,\varepsilon )$ the closed
ball of a radius $\varepsilon$ around $x$, i.e.
$$
B(x,\varepsilon)=\{ y \in \mathbb{R}^{n} \ | \ \| y-x \| \le \varepsilon \}.
$$

%==========================lm 2.1=============================

\begin{lemma} \label{lm:2.1}
Let  $\{p^i \}$ and $\{g^i \}$ be sequences in $\mathbb{R}^{n}$ such that
\begin{equation} \label{eq:2.1}
 p^{0} = g^{0}, \ p^{i+1} = {\rm Nr } \ {\rm conv} \{ p^{i}, g^{i+1} \}, \quad i=0,1,
                                                           \ldots
\end{equation}
Then
$$
p^{i} \in {\rm conv} \{ g^{0},\ldots,g^{i}\}, \quad i=0,1,\ldots
$$
\end{lemma}

%==========================lm 2.2=============================

\begin{lemma} \label{lm:2.2}
Let  $\{p^i \}$ and $\{g^i \}$ be sequences in $\mathbb{R}^{n}$ such that
(\ref{eq:2.1}) holds and
$$
  \| g^{i} \| \leq C < \infty,  \langle g^{i+1}, p^{i} \rangle \leq \theta \| p^{i}\|^{2},
  \quad i=0,1,\ldots , \ \theta \in (0, 1).
$$
Then
$$
 \|p^{i}\| \leq C/((1-\theta)\sqrt{i+1})   \quad {\rm for } \ i=0,1,\ldots
$$
\end{lemma}

We recall some concepts and properties from Convex Analysis.
If a function $f : \mathbb{R}^{n} \to \mathbb{R}$
is convex, one can define its subdifferential:
$$
 \partial f(x) = \{  g \in \mathbb{R}^{n}  \ \vrule \ f(y)-
  f(x) \geq  \langle g, y-x \rangle   \quad \forall y \in \mathbb{R}^{n} \},
$$
which is non-empty, convex and compact at any point $x$. Moreover,
the subdifferential mapping $x \mapsto \partial f(x)$ is upper semicontinuous
at any point $x$.

We need also a simple deviation estimate for convex nonsmooth functions; see
\cite[Lemma 4.3]{Kon13}.

%==========================lm 2.3=======================================

\begin{lemma} \label{lm:2.3}
Let
$$
p=\sum \limits _{j=1}^{m} \mu _{j} g^{j}, \ g^{j} \in \partial
f(y^{j}), \ y^{j} \in B(x, \delta ), \ \sum \limits_{j=1}^{m} \mu
_{j}=1, \ \mu _{j} \geq 0, \ j=1,\ldots, m.
$$
Then
\begin{equation} \label{eq:2.3}
\begin{array}{c}
\displaystyle f(y)-\sum \limits _{j=1}^{m} \mu _{j} f(y^{j}) \geq
 \langle p, y-x \rangle  -  \delta L_{x, \delta}
\end{array}
\end{equation}
for any $y \in \mathbb{R}^{n}$ where
$$
L_{x, \delta}= \max \limits_{j=1,\ldots, m} \|g^{j}\|.
$$
\end{lemma}
{\bf Proof.}  For any  $y $  we have
\begin{eqnarray*}
 f(y)-f(y^{j}) &\geq&  \langle g^{j}, y-y^{j} \rangle  \geq  \langle
            g^{j}, y-x \rangle  - \delta \|g^{j}\| \\
  &\geq&  \langle g^{j}, y-x \rangle  - \delta L_{x, \delta}
\end{eqnarray*}
for $j=1,\ldots, m$. Multiplying these inequalities by $\mu _{j}$ and summing
over $j=1,\ldots, m$, gives (\ref{eq:2.3}). \QED

%33333333333333333333333333333333333333333333333333333333333333333333333333

\section{The Basic Method and Its Convergence} \label{sc:3}

We will use the following set of basic assumptions.

{\bf (A1)}  {\em The function $f : \mathbb{R}^{n} \to \mathbb{R}$ is convex.}

{\bf (A2)} {\em There exists a number $\alpha $ such that the set
$$
E_{\alpha}=\left\{ x \in \mathbb{R}^{n} \ | \ f (x) \leq \alpha \right\}
$$
is non-empty and bounded.}

Clearly, {\bf (A2)} is a general coercivity condition, which implies that
the set $E_{\alpha}$ is bounded for any $\alpha$  if non-empty.
If {\bf (A1)} and {\bf (A2)} hold, problem (\ref{eq:1.1}) has a solution.

We now describe the non-monotone conjugate subgradient type method
for problem (\ref{eq:1.1}), which involves a
simple adaptive step-size procedure without line-search.

%===========================Method (CSGM)==================================
\medskip
\noindent {\bf Method (CSGM).}

{\em Step 0:} Choose a point $x^{0}$, numbers $\mu$ and  $\theta \in (0,1)$, and sequences
$\{\alpha _{s}\}$, $\{\beta _{m}\}$, $\{\eta _{t}\}$, $\{d _{t}\}$ such that
\begin{eqnarray}
  & &   \alpha _{s} \in (0,1), \{\alpha _{s}\} \to 0; \
    \beta _{m} > 0, \{\beta _{m}\} \to 0,  \sum \limits_{m=0}^{\infty } \beta _{m} = \infty ;  \label{eq:3.1}\\
   &&  \eta _{t} > 0, \{\eta_{t}\} \to 0; \ d _{t} > 0, \{d _{t}\} \to 0.
    \label{eq:3.2}
\end{eqnarray}
 Set $k:=0$, $l:=0$, $m:=0$, $s:=0$, $t:=0$, $b:=0$, $\lambda_{0} :=\beta _{0}$, $y^{0} := x^{0}$.
 Compute $ g^{0} \in \partial f (x^{0})$, set $ u^{0} :=x^{0}$, $ p^{0} :=g^{0}$.

{\em Step 1:}  If $\| p^{k}\| \le \eta_{t}$, set $ p^{k} :=g^{k}$,
 $l:= l+1$, $t:=t+1$,  $b:=0$. ({\em norm restart})

{\em Step 2:} Set $y^{k+1} := x^{k}-\lambda_{k}p^{k}$, $b := b+\lambda_{k}\| p^{k}\|$,
compute $ g^{k+1} \in \partial f (y^{k+1})$. If
\begin{equation} \label{eq:3.3}
  f(y^{k+1}) \leq f(x^{k})-\theta \lambda_{k}\|p^{k}\|^{2},
\end{equation}
set $x^{k+1}:=y^{k+1}$, $\lambda_{k+1} :=\lambda_{k}$
and go to Step 4. ({\em descent step})

{\em Step 3:}  Set $\lambda_{k+1} := \alpha _{s} \beta _{m}$, $s:=s+1$.
If $f (y^{k+1}) \leq \mu$, set $x^{k+1}:=y^{k+1}$ and go to Step 4.
Otherwise  take $ p^{k+1}:= g^{k+1}\in \partial f (u^{k})$,
 set $x^{k+1}:=u^{k}$, $u^{k+1}:=u^{k}$, $m:= m+1$, $\lambda_{k+1} := \beta _{m}$,
$k:=k+1$,  $t:=t+1$, $s:=0$,  $b:=0$  and go to Step 1. ({\em function value restart})

{\em Step 4:}  If $f(x^{k+1}) < f(u^{k})$, set $u^{k+1}:=x^{k+1}$.
If $b > d _{t}$, take
 $ p^{k+1}:= g^{k+1}\in \partial f (x^{k+1})$, set $m:= m+1$, $\lambda_{k+1} := \beta _{m}$,
$k:=k+1$,  $t:=t+1$, $s:=0$,   $b:=0$ and go to Step 1. ({\em distance restart})

{\em Step 5:} Set
$$
p^{k+1} := {\rm Nr} \ {\rm conv} \{ p^{k}, g^{k+1}\},
$$
$k:=k+1$ and go to Step 1.
\medskip

 According to the description, at each iteration,
 the current direction $p^{k}$ belongs to the convex hull of
 the subgradients from several previous iteration points.
 In the restart cases, we set $ p^{k}:= g^{k}\in \partial f (x^{k})$.
 The norm restart implies the decrease of the parameters
 $\eta_{t}$ and $d_{t}$ since the current point approximates a solution.
 The variable $l$ is a counter for the norm restarts.
 The distance and function value restarts imply the decrease of the parameters
 $\eta_{t}$ and $d_{t}$ and new starting step-size value $\beta _{m}$
 since the current approximation appeared not so precise.

 The strategy of descent steps and norm restarts follows the conjugate subgradient
 and adaptive gradient methods from \cite{Wol75,Kon84a,Kon18b}.
 The strategy of function value restarts follows the
  adaptive gradient method from \cite{Kon18b}.
 The  distance restart strategy follows in part the non monotone
 averaging subgradient methods from \cite[Chapter IV]{MGN87}.
Therefore, (CSGM) is an intermediate method since
it admits steps without descent, but descent and non descent steps imply
different choice of the next step-size values.
Observe that the sequence $\{u^{k}\}$ simply contains the best current points of
the sequence $\{x^{k}\}$, i.e.
$f(u^{k})=\varphi _{k}$ where $\varphi _{k} = \min \limits_{0 \leq i \leq k} f(x^{i})$.

In order to guarantee convergence of (CSGM)
we have to specialize the choice of the parameter
$\mu$. Take any $x^{*}\in X^{*}$, then there exists a ball
$B(x^{*},\rho)\supseteq E_{\alpha'}$ where $\alpha'=f(x^{0})$.
Let $\rho'$ denote the radius of the smallest ball
$B(x^{*},\rho')$  containing $E_{\alpha'}$ and let $E_{\mu'} \supset B(x^{*},\rho')$.

%3=========================thm 3.1============================================

\begin{theorem} \label{thm:3.1}
Let the assumptions {\bf (A1)}--{\bf (A2)} be fulfilled and
$\mu \geq \mu'$. Then:

(i) The sequence $\{x^k \}$ has a limit point, which belongs to the set
$X^{*}$.

(ii)  It holds that
\begin{equation}\label{eq:3.4}
\lim \limits_{k\rightarrow \infty} \varphi _{k}=f^{*}.
\end{equation}
\end{theorem}
{\bf Proof.}  First we observe that the sequence $\{x^{k}\}$ belongs to the bounded
set $E_{\mu}$ and must have limit points.

Let us consider several possible cases.

\textit{Case 1: The number of changes of the index $l$ is infinite.} \\
Then we have $t \to \infty$  and $\{d _{t}\} \to 0$ because of  (\ref{eq:3.2}).
Besides, in accordance with Lemma \ref{lm:2.1} there exist
subsequences of indices $\{k_{l}\}$ and $\{t_{l}\}$ such that
$$
 \|p^{k_{l}}\|\leq\eta_{l}, \
   p^{k_{l}} \in \mbox{conv}
    \bigcup_{\|y-x^{k_{l}}\| \leq d_{t_{l}} } \partial f(y), \ \{d _{t_{l}}\} \to 0,
$$
for $l=0,1,\ldots$ Take  an arbitrary limit point $x^{*}$  of the subsequence
$\{x^{k_{l}}\}$. Without loss of generality we suppose that
$$
\lim_{l \to \infty } x^{k_{l}}= x^{*}.
$$
Since the mapping $x \mapsto \partial f(x)$ is upper semicontinuous,
these properties together with (\ref{eq:3.2}) yield
$x^{*}\in X^{*}$. It follows that (\ref{eq:3.4}) is also true.

\textit{Case 2: The number of changes of the index $l$ is finite.} \\
This situation admits two possible cases.

\textit{Case 2a: The number of changes of the index $m$ is finite.} \\
Then  we have $t \leq t' < \infty$, $d _{t} \geq \bar d >0$  and
$\eta _{t} \geq \bar \eta >0$, hence
\begin{equation}\label{eq:3.5}
\|p^{k}\|\geq \bar \eta \quad  \mbox{for } \ l=0,1,\ldots
\end{equation}
If we suppose that the number of changes of the index $s$ is finite,
from (\ref{eq:3.3}) we have
$$
  f(x^{k+1}) \leq f(x^{k})-\theta \lambda_{k}(\bar \eta)^{2}
$$
for $k$ large enough. But now (\ref{eq:3.1}) gives
$$
  \sum \limits_{k=0}^{\infty } \lambda _{k} = \infty,
$$
hence $f(x^{k}) \to -\infty$ as $k \to \infty$, which is a contradiction.
Otherwise, the number of changes of the index $s$ is infinite.
Then Lemma \ref{lm:2.2} gives $ \|p^{k}\| \to 0$ as $k \to \infty$,
which contradicts (\ref{eq:3.5}).
That is, this case is impossible.

\textit{Case 2b: The number of changes of the index $m$ is infinite.} \\
Then we have $t \to \infty$  and $\{d _{t}\} \to 0$.
For a number $\varepsilon>0$ we set
$$
 T_{\varepsilon }=\left\{ {x \in \mathbb{R}^{n} \ | \
f(x) < f^{*}+\varepsilon } \right\}
$$
and
$$
B(T_{\varepsilon }, \varepsilon )=\{ {x \in \mathbb{R}^{n} \ | \
\inf_{u \in T_{\varepsilon }} \|x-u\| \leq {\varepsilon }   } \}.
$$
Let us suppose that there exists a number $\varepsilon>0$ such that
$x^{k} \notin B(T_{\varepsilon }, \varepsilon )$ for any $k=0, 1, \dots$
Take a point $x^{*}\in X^{*}$ related to the choice of $\rho'$, then
$$
\|x^{k+1}-x^{*}\|^{2}=\|x^{k}-x^{*}\|^{2}-2\lambda  _{k}
 \langle p^{k}, x^{k}-x^{*} \rangle + \lambda _{k}^{2} \|p^{k}\|^{2}.
$$
Next, at each point $x^{k}$ we have
$$
p^{k} \in {\rm conv} \{ g^{j}\}_{j\in J_k}, \ g^{j} \in \partial
f(y^{j}), \ y^{j} \in B(x^{k}, d _{t}),
$$
where $t$ is the current index value at the $k$-th iteration.
Take an index $k'$ large enough such that $d _{t} < \varepsilon $ if $k \geq k'$.
This means that  $y^{j} \notin T_{\varepsilon }$. Now from  Lemma \ref{lm:2.3}
we have
$$
 -\varepsilon  \geq  \langle p^{k}, x^{*}-x^{k} \rangle -
d _{t} C_{k},
$$
 where
$$
C_{k}= \max \limits_{j\in J_k} \|g^{j}\| \leq C < \infty.
$$
It follows that
\begin{equation}\label{eq:3.6}
\|x^{k+1}-x^{*}\|^{2}\leq \|x^{k}-x^{*}\|^{2}-\lambda  _{k} \varepsilon'
\end{equation}
for some $\varepsilon' \in (0, \varepsilon) $ and $k $ large enough.
Due the choice of the parameter $\mu$ this means that
all the points $x^{k}$ will be contained in the ball
$B(x^{*},\rho')$ and (\ref{eq:3.6}) holds for $k $ large enough. It follows that the
function value restart does not occur if $k \geq k''$ where $k'' \geq k'$ is chosen
large enough.
Since the number of changes of the index $m$ is infinite,
(\ref{eq:3.1}) gives
$$
  \sum \limits_{k=0}^{\infty } \lambda _{k} = \infty,
$$
hence $\|x^{k}-x^{*}\| <0$ for $k$ large enough, which is a contradiction.

Therefore, for each $\varepsilon>0$ there exists a number $k(\varepsilon)$ such that
$x^{k(\varepsilon)} \in B(T_{\varepsilon }, \varepsilon )$.
Taking any sequence $\{\varepsilon_{h}\} \to 0$, we obtain that any limit point of  the corresponding
sequence $\{x^{k(\varepsilon_{h})} \}$ belongs to $X^{*}$ and that
(\ref{eq:3.4}) holds true.
\QED

%4444444444444444444444444444444444444444444444444444444444444444444444444444444444444444444444444

\section{Modifications and Implementation}\label{sc:4}

The above descent method admits various modifications and
extensions. For instance, we can determine the sequence $\{\eta _{l}\}$ to be dependent
only of the index $l$ and take
$$
\eta _{l} > 0, \{\eta_{l}\} \to 0
$$
instead of (\ref{eq:3.2}). Besides, in Step 1 we should
check now the inequality $\| p^{k}\| \le \eta_{l}$.
Then all the assertions of Theorem \ref{thm:3.1}
remain true with small modifications in the proof.
Nevertheless, we think that the previous version is more suitable for application
since it combines rules from descent and non descent subgradient methods for all the
parameters. Next, it is not necessary in fact to evaluate $\mu'$ directly.
We can first take any $\mu \geq f(x^{0})$ and fix a number $\alpha''>0$.
Afterwards, we should set $\mu:=\mu +\alpha''$ if
the function value restart occurs.

We now illustrate opportunities of (CSGM) with its behavior
in the smooth case.
Let us suppose that the gradient of $f$ is Lipschitz continuous with constant $L$.
We first fix temporarily the values $\eta _{t}= +\infty$ and $d _{t}= +\infty$.
Then (CSGM) has norm restarts at each iteration and behaves as
the adaptive gradient method from \cite{Kon18b}. Next, then
by (\ref{eq:3.3}) we have
$$
  f(x^{k+1}) \leq f(x^{k})-\theta \lambda_{k}\|f'(x^{k})\|^{2}
$$
for $k$ large enough since this inequality holds
if $ \lambda_{k}\leq \bar  \lambda =2(1-\theta)/L$; see \cite[Chapter II, \S 1]{PD78}.
This means that $ \lambda_{k}=  \lambda'>0$ for $k$ large enough.
If we suppose in addition that  $f$ is strongly convex with constant $\varkappa$, then
$$
 f(x^{k}) -f^{*} \leq (2/\varkappa) \|f'(x^{k})\|^{2},
$$
and the above method converges linearly to the optimal value $f^{*}$;
see \cite[Chapter II, \S 1]{PD78}.
Therefore, all these properties will hold if we choose
$\eta _{0}$ and $d _{0}$ large enough and take these parameters
in conformity with (\ref{eq:3.2}) but tending to $0$ slower than linearly.
It follows that   divergent series rules for parameters do not prevent in general from
a linear rate of convergence  if they involve adaptive step-size procedures.
 At the same time, averaging gradients may enhance performance
after proper norm restart regulation with evaluation of the goal function properties along
current points. For this reason, we should avoid the situation where
$d _{t} \approx \beta _{m}\| p^{k}\|$ in the method
in order to guarantee averaging gradient iterations.

We now intend to somewhat specialize the choice of
parameters in (CSGM) and describe its implementation, which involves two different rules after
norm and distance restarts. That is, the rule after a norm restart
corresponds to those in descent subgradient and adaptive gradient methods, whereas
the rule after a distance restart follows that in non descent subgradient methods.

\newpage

%===========================Method (CSGI)==================================
\medskip
\noindent {\bf Method (CSGI).}

{\em Step 0:} Choose a point $x^{0}$, numbers $\mu$ and  $\theta \in (0,1)$, and sequences
$\{\alpha' _{s}\}$, $\{\alpha'' _{l}\}$, $\{\beta' _{m}\}$, $\{\beta'' _{m}\}$, $\{\beta''' _{m}\}$ such that
\begin{eqnarray}
  & &   \alpha' _{s} \in (0,1), \{\alpha' _{s}\} \to 0; \  \alpha'' _{l} \in (0,1), \{\alpha'' _{l}\} \to 0; \label{eq:4.1}\\
  & &   \beta' _{m} > 0, \{\beta' _{m}\} \to 0,  \sum \limits_{m=0}^{\infty } \beta' _{m} = \infty ; \
   \beta'' _{m} > 0, \{\beta'' _{m}\} \to 0;  \label{eq:4.2}\\
   && \beta''' _{m} > 0, \{\beta''' _{m}\} \to 0.
    \label{eq:4.3}
\end{eqnarray}
 Set $k:=0$, $l:=0$, $m:=0$, $s:=0$, $t:=0$, $b:=0$, $\lambda_{0} :=\beta' _{0}$,
 $\eta_{0} :=\beta'' _{0}$,$d_{0} :=\beta''' _{0}$, $y^{0} := x^{0}$.
 Compute $ g^{0} \in \partial f (x^{0})$, set $ u^{0} :=x^{0}$, $ p^{0} :=g^{0}$.

{\em Step 1:}  If $\| p^{k}\| \le \eta_{t}$, set $ p^{k} :=g^{k}$, $\eta_{t+1} := \alpha'' _{l} \beta'' _{m}$,
$d_{t+1} := \alpha'' _{l} \beta''' _{m}$,
 $l:= l+1$, $t:=t+1$,  $b:=0$. ({\em norm restart})

{\em Step 2:} Set $y^{k+1} := x^{k}-\lambda_{k}p^{k}$, $b := b+\lambda_{k}\| p^{k}\|$,
compute $ g^{k+1} \in \partial f (y^{k+1})$. If
$$
  f(y^{k+1}) \leq f(x^{k})-\theta \lambda_{k}\|p^{k}\|^{2},
$$
set $x^{k+1}:=y^{k+1}$, $\lambda_{k+1} :=\lambda_{k}$
and go to Step 4. ({\em descent step})

{\em Step 3:}  Set $\lambda_{k+1} := \alpha' _{s} \beta _{m}$, $s:=s+1$.
If $f (y^{k+1}) \leq \mu$, set $x^{k+1}:=y^{k+1}$ and go to Step 4.
Otherwise  take $ p^{k+1}:= g^{k+1}\in \partial f (u^{k})$,
 set $x^{k+1}:=u^{k}$, $u^{k+1}:=u^{k}$, $m:= m+1$, $\lambda_{k+1} := \beta' _{m}$,
 $\eta_{t+1} :=\beta'' _{m}$, $d_{t+1} :=\beta''' _{m}$,
$k:=k+1$,  $t:=t+1$, $s:=0$, $l:=0$,  $b:=0$  and go to Step 1. ({\em function value restart})

{\em Step 4:} If $f(x^{k+1}) < f(u^{k})$, set $u^{k+1}:=x^{k+1}$.  If $b > d _{t}$, take
 $ p^{k+1}:= g^{k+1}\in \partial f (x^{k+1})$, set $m:= m+1$, $\lambda_{k+1} := \beta' _{m}$,
 $\eta_{t+1} :=\beta'' _{m}$,$d_{t+1} :=\beta''' _{m}$,
 $k:=k+1$, $t:=t+1$, $s:=0$, $l:=0$,   $b:=0$ and go to Step 1. ({\em distance restart})

{\em Step 5:} Set
$$
p^{k+1} := {\rm Nr} \ {\rm conv} \{ p^{k}, g^{k+1}\},
$$
$k:=k+1$ and go to Step 1.
\medskip

Following the lines of the proof  of Theorem \ref{thm:3.1},
we can conclude that the same assertions are true for (CSGI).

%4=========================thm 4.1============================================

\begin{theorem} \label{thm:4.1}
Let the assumptions {\bf (A1)}--{\bf (A2)} be fulfilled and
$\mu \geq \mu'$. Then
the sequence $\{x^k \}$ generated by (CSGI) has a limit point, which belongs to the set
$X^{*}$, besides, relation (\ref{eq:3.4}) holds.
\end{theorem}

Although (CSGI) involves several sequences of
parameters, it seems reasonable to apply the same rule for $\{\alpha' _{s}\}$ and $\{\alpha'' _{l}\}$,
as well as the same rule for $\{\beta' _{m}\}$, $\{\beta'' _{m}\}$ and $\{\beta''' _{m}\}$.
This in particular means that
$$
 \sum \limits_{m=0}^{\infty } \beta'' _{m} = \infty , \ \sum \limits_{m=0}^{\infty } \beta''' _{m} = \infty .
$$
However, these additional rules are not obligatory for providing convergence.

%55555555555555555555555555555555555555555555555555555555555555555555555555555555555555555

\section{Computational Experiments}\label{sc:5}

In order to check the performance of the proposed method we carried
out computational experiments.  The main goal was to compare it
with subgradient methods without line-search, which  have
similar storage and computation expenses per iteration.

We first took the simplest subgradient method (\ref{eq:1.2})--(\ref{eq:1.3}) abbreviated as
(SGM), the non monotone averaging subgradient method from \cite{Che81}
abbreviated as (NASGM), the simple descent subgradient method
from \cite{Kon82,Kon84a} abbreviated as (DSGM),
and the proposed method (CSGI).
We took for comparison the well-known non-smooth
test problem from \cite{ShS72} (see also \cite{Sho85,Pol87})
with $n=5$, the optimal value $f^{*}=22.60016$,
and the starting point $x^{0} = (0,0,0,0,1)^{\top} $.

We compared all the methods for different accuracy $\varepsilon$
with respect to the goal function deviation
$$
\Delta (x)=f(x)-f^{*}.
$$
They were implemented in Delphi with double precision
arithmetic. Namely, we indicate
the total number of iterations (it)
(or the total number of subgradient calculations) for attaining
the desired accuracy $\varepsilon$.

We took the rule $\lambda_{k} =\lambda/(k+1)$ with $\lambda =0.1$ for  (SGM).
Since (NASGM) involves several sequences of parameters tending to $0$,
we decreased them linearly with the same ratio $\sigma = 0.9$.
For (DSGM), we also used the same rule with the ratio $\sigma = 0.6$
for two sequences of parameters tending to $0$. For
(CSGI), we set $\theta =0.3$  and the rules $\beta' _{m} =\beta' /(m+1)$,
$\beta'' _{m} =\beta'' /(m+1)$, $\beta''' _{m} =\beta''' /(m+1)$
with $\beta'=0.05$, $\beta''=0.4 \| g^{0}\|$, and $\beta'''=\beta' \| g^{0}\|/0.7$.
Also, we set $\alpha' _{0}=\alpha'' _{0}=\sigma=0.8$ and the rule
$$
\alpha' _{s+1} =\sigma \alpha' _{s}, \ \alpha'' _{l+1} =\sigma \alpha'' _{l}.
$$
For simplicity, we implemented (CSGI) with $\mu=+\infty$, i.e. without 
function value restarts.  The results  are given in Table \ref{tbl:1}.

\begin{table} \caption{Comparison of subgradient methods} \label{tbl:1}
\centering
\begin{tabular}[t]{|r|r|r|r|r|r|r|r|}
\hline
\multicolumn{2}{|c|}{SGM} & \multicolumn{2}{|c|}{NASGM} & \multicolumn{2}{|c|}{DSGM} & \multicolumn{2}{|c|}{CSGI}\\
\hline
   {$\varepsilon$} & {it} & {$\varepsilon$} & {it} & {$\varepsilon$} & {it} & {$\varepsilon$} & {it}\\					
\hline
0.1 & 81 & 0.1 & 30 & 0.1 & 92 & 0.1 & 141\\					
\hline
 0.01 & 320 & 0.01 & 63 & 0.01 & 352 & 0.01 & 253\\
\hline
 0.001 & 1645 & 0.004 & 10000 & 0.001 & 1058 & 0.001 & 466\\	
\hline
 0.0001 & 8243  & -    & - & 0.0001 & 2809 & 0.0001 & 640\\	
\hline
 0.00002 & 35000 & -   & - & 0.00001 & 5909 & 0.00001 & 860\\	
\hline
\end{tabular}
\end{table}

The implementation of (CSGI) showed rather rapid convergence
in comparison with the other methods.
Convergence of (NASGM) appeared rather instable. At the same time,
convergence of (DSGM) appeared not so rapid, but stable.

Somewhat different strategies for determining step-sizes in subgradient methods
with averaging procedures were proposed in \cite{Nes09,NeS15}. They were deduced from
theoretical convergence rates for these methods on the class of convex minimization
problems. Besides, comparative computational experiments with three methods were also described
in \cite{NeS15}. They involved  the same subgradient method (\ref{eq:1.2})--(\ref{eq:1.3})
with the best theoretical rule:
\begin{equation} \label{eq:5.1}
  \lambda _{k}:=\lambda/\sqrt{(k+1)}, \quad \lambda =\|x^{0}-x^{*}\|/L,
\end{equation}
where $L$ is an upper bound for the norm of subgradients, $x^{*}\in X^{*}$. In what follows,
we abbreviate this method as (SGMT). Next, the so-called
 simple dual averaging method from \cite{Nes09} was taken. It can be written as follows:
$$
 x^{k+1} :=x^{0} - \lambda _{k} p^{k}, \  p^{k}:=\sum \limits_{k=0}^{\infty } g^{k}, \
 g^{k} \in \partial f (x^{k}),
$$
where $ \lambda _{k}$ was chosen as in (\ref{eq:5.1}).
We abbreviate this method as (ASG). In addition, the so-called
 method of simple double averaging from \cite{NeS15} was used. It  can be written as follows:
\begin{eqnarray*}
 & &  x^{k+1} :=\mu _{k}x^{k}+ (1 -\mu _{k}) y^{k}, \ \mu _{k}:=(k+1)/(k+2),     \\
  & &   y^{k} :=x^{0} - \lambda _{k} p^{k}, \ p^{k}:=\sum \limits_{k=0}^{\infty } g^{k}, \ g^{k} \in \partial f (x^{k}),
\end{eqnarray*}
where $ \lambda _{k}$ was chosen as in (\ref{eq:5.1}).
We abbreviate this method as (DASG).

We wrote programs for these three methods and
took for comparison the above test problem with
the same starting point $x^{0} = (0,0,0,0,1)^{\top} $.
Also, we took $L= \| g^{0}\|$. The results  are given in Table \ref{tbl:2}.

\begin{table} \caption{Comparison of subgradient methods} \label{tbl:2}
\centering
\begin{tabular}[t]{|r|r|r|r|r|r|}
\hline
\multicolumn{2}{|c|}{SGMT} & \multicolumn{2}{|c|}{ASG} & \multicolumn{2}{|c|}{DASG} \\
\hline
   {$\varepsilon$} & {it} & {$\varepsilon$} & {it} & {$\varepsilon$} & {it} \\					
\hline
0.1 & 116 & 2.038 & 10000 & 0.1 & 324 \\					
\hline
 0.01 & 4510 & - & - & 0.01 & 3254 \\
\hline
 0.0013 & 35000 & - & - & 0.001 & 34169 \\	
\hline
\end{tabular}
\end{table}

From the calculations we can conclude that (SGMT) is slow
in comparison with (SGM), but (ASG) appeared very slow.
At the same time, convergence of (DASG) appeared better than (SGMT),
but was not so rapid even in comparison with (SGM).
Therefore, these variants of methods with averaging subgradients
demonstrated rather slow convergence. In fact, their step-sizes
strategies stem from the worst case analysis of convergence
for the whole class of minimization problems. At the same time, any
nonlinear function may behave in a different manner on
different sets containing iteration points.
Hence, an iterative solution method should utilize
some adaptive parameter strategies for attaining
better convergence properties.

%%%%%%%%%%%%%%%%%%%%%%%%%%%%%%%%%%%%%%%%%%%%%%%%%%%%%%%%%%%%%%%%%%%%%%%%%%%%%%%%%%%%%

\section*{Acknowledgement}

The results of this work were obtained within the state assignment of the
Ministry of Science and Education of Russia, project No. 1.460.2016/1.4.
This work was supported by Russian Foundation for Basic Research, project No.
19-01-00431.

%##########################################################################

\end{document}